\documentclass[a4paper,leqno,11pt]{article}
\usepackage{latexsym}
\usepackage{amsthm,amsmath,amssymb}

\usepackage[dvips]{graphicx, color}

\newtheorem{thm}{Theorem}
\newtheorem{mthm}{Theorem}
\newtheorem{mdef}{Definition}
\newtheorem{mmthm}{Theorem$'$}
\newtheorem{prop}[thm]{Proposition}

\newtheorem{axiom}[]{Axiom}

\begin{document}
\title{\bf {\Large{Geometry of Numbers}}}
\author{\bf Lin WENG}
\date{\bf }
\maketitle

\noindent
{\bf Abstract.} {\footnotesize We develop a global cohomology theory for number fields by offering
topological cohomology groups,
an arithmetical duality, a Riemann-Roch  theorem, 
and two types of vanishing theorem. As applications, we study moduli spaces of
semi-stable lattices, and introduce non-abelian zeta functions for number fields.\footnote{
 Parts of this paper were circulated under the title: {\it Riemann-Roch, Stability and
New Non-Abelian Zeta Functions for Number Fields}, arXiv:math/0007146. The new developments on vanishing theorem and strong stability  are added to complete the theory.}}
\vskip 1.0cm
\section{\bf  Global Cohomology}
\subsection{\bf Adelic cohomology groups} 
Let $F$ be a number field, and $\mathcal O_F$ its integer ring. Denote $S:=S_{\mathrm{fin}}\cup S_\infty$ 
the set of normalized non-Archimedean and Archimedean places of $F$. 
For all $v\in S$, write $F_v$ the $v$-completion of $F$, $\mathcal O_v$ its integer ring. 

Denote $\mathbb A$  the adelic ring of $F$, $GL_n(\mathbb A)$ the associated general linear group.
For each $g=(g_\frak p;g_v)\in GL_n(\mathbb A)$, $\frak p\in S_{\mathrm{fin}}, v\in S_\infty$,
 introduce an auxiliary topological space $\mathbb A^n(g)$ as follows:

\noindent
(i) Set theoretically,
$$\mathbb A^n(g):=\Bigg\{{\bf x}\in\mathbb A^n\,\Big|\,\begin{matrix}\exists {\bf a}\in F^n, \ \mathrm{s.t.}\ {\bf x_v}={\bf a}, &\forall v\in S_\infty\\
g_\frak p\cdot{\bf  x}_\frak p\in\mathcal O_{\frak p}^n,\,g_\frak p\cdot {\bf a}\in\mathcal O_{\frak p}^n,
&\forall \frak p\in S_{\mathrm{fin}}\end{matrix}\Bigg\};$$ 

\noindent
(ii) Topologically, first introduce a new topological structure on $\mathbb A^n$ by keeping the finite part but altering
 its metric at $v\in S_\infty$ using the positive definite matrix $g_\sigma^t\cdot g_\sigma$ (resp. $\bar g^{t}_\tau\cdot g_\tau$) when $v=\sigma$ is real (resp. $v=\tau$ is complex); 
 then equip $\mathbb A^n(g)$ with the induced topological structure from
the embedding $\mathbb A^n(g)\subset\mathbb A^n$. 

\begin{mdef} For a matrix idele $g\in GL_n(\mathbb A)$, define
the {\it 0-th} and the {\it 1-st arithmetical
cohomology groups} of $g$ by:

\noindent
(1) As abstract abelian groups,
  $$\displaystyle{H^0(F,g):=\mathbb A^n(g)\cap K^n,\qquad\mathrm{and}\qquad
H^1(F,g):=\mathbb A^n\big/\big(\mathbb A^n(g)+ K^n\big);}$$

\noindent
(2) As topological spaces, $H^0(F,g)$ and $H^1(F,g)$ are equipped with topologies induced from the altered $\mathbb A^n$. 
\end{mdef}
\eject
\begin{prop} As locally compact topological spaces, $H^0(F,g)$ is discrete and $H^1(F,g)$ is compact.
\end{prop}

\noindent
{\it Proof.} $\bullet$ {\bf Discreteness of} $H^0$: Following Minkowski, via Archimedean places, embed $F$ into  $\mathbb R^{r_1}\times\mathbb C^{r_2}$, where $r_1$ (resp. $r_2$) denotes the number of real places (resp. complex places).
 Introduce 
$$H^0_{\mathrm{fin}}(F,g):=\Big\{{\bf x}\in F^n\,\big|\,
g_\frak p\cdot{\bf  x}\in\mathcal O_{\frak p}^n, \forall \frak p\in S_{\mathrm{fin}}\Big\}.$$
By definition, particularly, by paying attention to the infinite components, we conclude that,
 as abelian groups, $H^0(F,g)\simeq H^0_{\mathrm{fin}}(F,g),$ and the natural embedding
$H^0(F,g)\hookrightarrow \mathbb A_\infty^n,$ via $$H^0(F,g)\simeq H^0_{\mathrm{fin}}(F,g)\hookrightarrow F^n\buildrel i\over\hookrightarrow 
\Big(\mathbb R^{r_1}\times\mathbb C^{r_2}\Big)^n=\mathbb A_\infty^n,$$ induces a natural discrete topological structure on $H^0(F,g)$.

\noindent
$\bullet$ {\bf Compactness of} $H^1$: This is equivalent to the {\it strong approximation theorem} for adeles. 
Indeed, if we equip 
$\mathbb A^n_\infty\simeq\prod_{\sigma:\mathbb R}\mathbb R^n\times\prod_{\tau:\mathbb R}\mathbb C^n$ with the twisted metric structure $\rho(g)$ obtained from the
positive definite matrix $(g^t_\sigma\cdot g_\sigma;\overline g_\tau^t\cdot g_\tau)$, then
$H^0(F,g)\simeq H^0_{\mathrm{fin}}(F,g)$ can be viewed as a full rank lattice $\Lambda(g)$ within.

\begin{prop} As locally compact topological groups,
we have the following natural isomorphism $$H^1(F,g)\simeq \Big(\prod_{\sigma:\mathbb R}\mathbb R^n\times\prod_{\tau:\mathbb R}\mathbb C^n,\rho(g_\infty)\Big)\Big/ \Lambda(g).$$
\end{prop}

\noindent
{\it Proof.} By keeping track on the metrics, it suffices to show that, as abstract groups, we have an isomorphism 
$$\mathbb A^n/(\mathbb A^n(g)+F^n)\,\simeq\, \mathbb A_\infty^n/i(H^0_{\mathrm{fin}}(F,g)).$$
For this, we use the natural morphism
$\phi:\mathbb A_\infty^n\to \mathbb A^n/(\mathbb A^n(g)+F^n)$ induced from
the embedding $\mathbb A_\infty\hookrightarrow \mathbb A$.  Hence we need to prove that

\hskip 1.50cm{\it (i)\ $\phi$\ is\ surjective;\ $\qquad\&\qquad$\ (ii)\ $\mathrm{Ker}\,\phi=i(H^0_{\mathrm{fin}}(F,g))$.}

The surjectivity (i) is equivalent to  $\mathbb A^n_{\mathrm{fin}}\subset \mathbb A^n(g)+F^n$. This can be established using
the strong approximation property: Indeed, for any ${ g}=(g_\frak p;g_v)\in \mathbb A^n$, 
there exists a finite set $S_0\subset S_{\mathrm{fin}}$ such that $g_\frak q\in GL_n(\mathcal O_\frak q)$ for $\frak q\not\in S_0$.
Moreover, by the strong approximation, there exists  ${\bf a}\in K^n$ such that

(a) $g_\frak p\cdot a_\frak p-{\bf a}\in\mathcal O^n_\frak p, \forall \frak p\in S_0$; and

(b) ${\bf a}\in \mathcal O_\frak q^n,\forall \frak q\not\in S_0$. 

(ii) is obtained by examining the definition of $\mathbb A^n(g)$ carefully. Indeed, this is
a direct consequence of the conditions for all components, particularly, those  at infinity, of the elements
in $H^0(F,g)$.

\subsection{Nine diagram: Riemann-Roch axioms}

Fix a finite place $\frak p$ of $F$. For $g\in GL_n(F_\frak p)$, define its {\it 0-th and 1-st cohomology group} by
 $$H^0(F_\frak p,g):=\{x\in F^n_\frak p: g\cdot x\in\mathcal O_\frak p\}\quad\mathrm{and}\quad  H^1(F_\frak p,g):=F^n_\frak p/H^0(F_\frak p,g)+F_\frak p^n.$$ Elements $g_1, g_2$ of $GL_n(F_\frak p)$ are called {\it equivalent}, denoted by $g_1\sim g_2,$  if $H^0(F_\frak p,g_1)=H^0(F_\frak p,g_2)$.
Clearly, $g_1\sim g_2$ if and only if there exists $g\in GL_n(\mathcal O_\frak p)$ such that $g_2=g\cdot g_1$.
Introduce then the quotient set $GL_n(F_\frak p)/GL_n(\mathcal O_\frak p)$ and a partial order $\leq$ on 
$GL_n(\mathcal O_\frak p)\backslash GL_n(F_\frak p)$: $$[g_1]\leq [g_2]\qquad\mathrm{if}\qquad
H^0(F_\frak p,g_1)\subset H^0(F_\frak p,g_2).$$

Globalizing this, we then obtain an equivalence relation $\sim$ on $GL_n(\mathbb A_\mathrm{fin})$ (and more generally on
$GL_n(\mathbb A_\mathrm{fin})\times\{g_\infty\}$ for a fixed  $g_\infty\in GL_n(\mathbb A_\infty)$)
and a partial order $\leq$ on
 $GL_n(\prod_{\frak p\in S_{\mathrm{fin}}}\mathcal O_\frak p)\backslash GL_n(\mathbb A_{\mathrm{fin}})$
such that $g_1\simeq g_2$ iff $H^0_{\mathrm{fin}}(F,g_1)=H^0_{\mathrm{fin}}(F,g_2)$ iff
there exists $g\in GL_n(\prod_{\frak p\in S_{\mathrm{fin}}}\mathcal O_\frak p)$ such that $g_2=g\cdot g_1$; and
$$[g_1]\leq [g_2]\qquad\mathrm{ iff}\qquad H^0_{\mathrm{fin}}(F,g_1)\subset H^0_{\mathrm{fin}}(F,g_2).$$

Recall that for $g\in GL_n(\mathbb A)$, we have the following 9-diagram with rows and columns exact:
$$\begin{matrix}
&&0&&0&&0&&\\
&&\downarrow&&\downarrow&&\downarrow&&\\
0&\to&H^0(F,g)&\to&\mathbb A^n(g)&\to& \mathbb A^n(g)/\mathbb A^n(g)\cap F^n&\to&0\\
&&\downarrow&&\downarrow&&\downarrow&&\\
0&\to&F^n&\to&\mathbb A^n&\to& \mathbb A^n/ F^n&\to&0\\
&&\downarrow&&\downarrow&&\downarrow&&\\
0&\to&F^n/\mathbb A^n(g)\cap F^n&\to&\mathbb A^n/\mathbb A^n(g)&\to &H^1(F,g)&\to&0\\
&&\downarrow&&\downarrow&&\downarrow&&\\
&&0&&0&&0.&&
\end{matrix}$$

Consequently, if $g=(g_\frak p;g_v), g'=(g_\frak p';g_v')$ satisfying $g_v=g_v', \forall v\in S_\infty$ and $[(g_\frak p)]\leq [(g'_\frak p)]$,  we have the following exact sequences by noticing that in the above diagram, the middle row remains invariant:
$$
H^0(F,g)\hookrightarrow H^0(F,g')\to
\mathbb A^n(g')/\mathbb A^n(g)\to H^1(F,g)\twoheadrightarrow H^1(F,g').$$
Note that with $g_\infty$ fixed, $\mathbb A^n(g')/\mathbb A^n(g)$ measures only finite places contributions. Namely, 
$\mathrm{deg}_{\mathrm{fin}}(g')-\mathrm{deg}_{\mathrm{fin}}(g')=\mathrm{deg}(g')-\mathrm{deg}(g')$ where 
$$\mathrm{deg}_{\mathrm{fin}}(g):=\mathrm{deg}_{\mathrm{fin}}({\mathrm{det}\,g)=\sum_\frak p
\mathrm{ord}_\frak p(\mathrm{det}g_\frak p})\cdot N_\frak p\log p.$$ 
This then further suggests a way to count $H^0$ and $H^1$. Namely, it should bear the rules of

\begin{axiom} (Weak Riemann-Roch for  a fixed $g_\infty$) For $g$, $g'$ with a fixed $g_\infty$,
$$h^0(F,g')-h^1(F,g')-\mathrm{deg}(g')=h^0(F,g)-h^1(F,g)-\mathrm{deg}(g).$$
\end{axiom}

\begin{axiom} (Riemann-Roch for a fixed $n$) \newline
$ \hskip 2.0cm h^0(F,g)-h^1(F,g)-\mathrm{deg}(g)$ is independent of $g$.
\end{axiom}

Assuming it, then we may write this invariant of $F$ as $\Delta(F,n).$ By looking at diagonal elements, it is natural to introduce
the following
\begin{axiom} (Strong Riemann-Roch)
$$\Delta(F,n)=n\cdot\frac{1}{2}\Delta(F).$$
\end{axiom}

\subsection{Duality: Local and Global Pairings}
  
To establish a canonical duality between $H^0$ and $H^1$, we start with the natural pairing $$\langle\cdot,\cdot\rangle:\mathbb A\times \mathbb A\to \mathbb R/\mathbb Z\simeq \mathbb S^1\subset\mathbb C^*$$ defined by
$$\langle{\bf x},{\bf y}\rangle:=\sum_{v\in S}\langle x_v,y_v\rangle_v,\qquad \forall\ {\bf x}=(x_v),\, {\bf y}=(y_v)\in\mathbb A.$$
Here
$$\langle x_v,y_v\rangle_v:=\prod_{i=1}^ne^{2\pi \sqrt{-1}\cdot\chi_v(x_{v,i}\cdot y_{v,i})}$$ with $\chi_v:=\lambda_v\circ\mathrm{Tr}_{F_v/\mathbb Q_v}$, and $$\lambda_v:=\begin{cases}\mathbb Q_v\twoheadrightarrow
\mathbb Q_v/\mathbb Z_v\hookrightarrow \mathbb Q/\mathbb Z\hookrightarrow \mathbb R/\mathbb Z,&v=\frak p\in S_{\mathrm{fin}}\\
\mathbb R\twoheadrightarrow \mathbb R/\mathbb Z,&v\in F_\infty.\end{cases}$$
It is well-known that the above global pairing induces a natural isomorphism $\widehat{\mathbb A}^n\simeq \mathbb A^n$ (as locally compact groups) which in particular induces an isomorphism $(F^n)^\perp\simeq F^n$
(as discrete subgroups). On the other hand, by a direct local calculation from the definition ([Ta]), we conclude that, for local pairings 
$\langle\cdot,\cdot\rangle_\frak p$ of finite places $\frak p$,  $\frak a_\frak p^\perp\simeq\frak a^{-1}_\frak p\cdot\partial_\frak p$,
 where $\partial_\frak p$ denotes the local differential module at $\frak p$ (dual to $\mathcal O_\frak p$). 
For a fixed uniformizer $\pi_\frak p$  of $\mathcal O_\frak p$, denote  $\kappa_F$  the idele 
$(\pi_\frak p^{-\mathrm{ord}_\frak p(\partial_\frak p)};1)$. Since $H^0$ is insensitive towards local and global units,
the cohomology group $H^0(F,\kappa_F\cdot g^{-1})$ is well-defined.
Here, as usual, $ g^{-1}$ denotes the  inverse of  $g$.
Consequently, we have the following

\begin{prop} (Topological Duality) The global pairing $\langle\cdot,\cdot\rangle:\mathbb A\times \mathbb A\to \mathbb C^*$ 
induces a natural isomorphism between locally compact groups
$$\widehat{H^1(F,g)}\simeq H^0(F,\kappa_F\otimes g^{-1}).$$
\end{prop}

In fact, the above exposes a much stronger relation between $H^i$.
To understand this properly, let us recall that by Prop. 2, we have a refined isomorphism
$$H^1(F,g)\simeq \Big(\prod_{\sigma:\mathbb R}\mathbb R^n\times\prod_{\tau:\mathbb R}\mathbb C^n,\rho(g_\infty)\Big)\Big/ \Lambda(g).$$ 
Moreover, the above proof shows that $\Lambda(g)^\perp\simeq \Lambda(\kappa_F\cdot  g^{-1}).$
This completes the proof of the following:

\begin{mthm} (Arithmetic Duality) There is an isometry of $\mathcal O_F$-lattices
$$\widehat{H^1(F,g)}\simeq H^0(F,\kappa_F\cdot  g^{-1}).$$
\end{mthm}

\subsection{Riemann-Roch Theorem}
With the above duality, and the axioms for Riemann-Roch in mind, to have an Riemann-Roch for our setting, we need to introduce arithmetic
 counts $h^0$ and $h^1$ for discrete groups $H^0$ and compact groups $H^1$. Now the Riemann-Roch in Arakelov theory ([La1,2]) claims that $$\mathrm{deg}(g)-\frac{n}{2}\log\Delta_F=\chi(F,g):=-\log\Big(\mathrm{Vol}\Big(\prod_{\sigma:\mathbb R}\mathbb R^n\times\prod_{\tau:\mathbb R}\mathbb C^n,\rho(g_\infty)\Big)\Big/ \Lambda(g)\Big).$$ It is compatible with all axioms above.
Therefore, $h^i$ to be introduced should satisfy
\vskip 0.20cm
\noindent
({\bf 1}) a numerical duality
$$h^1(F,g)=h^0(F,\kappa_F\cdot  g^{-1});$$ 

\noindent
({\bf 2})  $h^0(F,g)-h^1(F,g)=-\log\mathrm{Vol}\Big(\prod_{\sigma:\mathbb R}\mathbb R^n\times\prod_{\tau:\mathbb R}\mathbb C^n,\rho(g_\infty)\Big)\Big/ \Lambda(g).$
\vskip 0.20cm
\noindent
$\bullet\ h^0(F,g)$: Count  $H^0(F,g)$ with the weight function
$${\bf e}_{\mathbb A}:=\prod_{\frak p\in S_{\mathrm{fin}}}{\bf 1}_{\mathcal O_\frak p^n}\times\prod_{v\in S_\infty}e^{-\pi\cdot N_v \cdot\|*\|_{\mathrm{can},v}}.$$ Here $N_v:= [F_v:\mathbb R]$ and
$\|*\|_{\mathrm{can},v}$ denotes the canonical metric at the place $v$. That is to say, define the 
{\it arithmetical count} of the discrete $H^0(F,g)$ by
$$\begin{aligned}\#_{\mathrm{Ar}}(H^0(F,g)):=&\int_{{\bf x}\in H^0(F,g)}{\bf e}_{\mathbb A}(g\cdot {\bf x})\,d\mu({\bf x})\\
=&
\sum_{{\bf x}\in H^0(F,g)}\prod_{\frak p\in S_{\mathrm{fin}}}{\bf 1}_{\mathcal O_\frak p^n}(g_\frak p\cdot x_\frak p)\times
\prod_{v\in S_\infty}e^{-\pi N_v\|g_v\cdot x_v\|_{\mathrm{can},v}}.\end{aligned}$$  
(Note that a shift by the multiplicative factor $g$ is taking place here.) 

\begin{mdef} Define the 
 0-th numerical  cohomology $h^0(F,g)$ of a matrix idele $g\in GL_n(\mathbb A)$ by 
$$\begin{aligned}h^0(F,G):=&\log \Big(\#_{\mathrm{Ar}}(H^0(F,g))\Big)=\log\Big(\sum_{{\bf x}\in H^0(F,g)}{\bf e}_{\mathbb A}(g\cdot {\bf x})\Big)\\
=&\log\Big(\sum_{{\bf x}\in H^0(F,g)}\prod_{\frak p\in S_{\mathrm{fin}}}{\bf 1}_{\mathcal O_\frak p^n}(g_\frak p\cdot x_\frak p)\times
\prod_{v\in S_\infty}e^{-\pi N_v \|g_v\cdot x_v\|_{\mathrm{can},v}}\Big).\end{aligned}$$
\end{mdef}

\noindent
$\bullet\ h^1(F,g)$: We start with the following natural

\begin{axiom} If $G$ is a discrete or compact group, then arithmetic counts $\#_{\mathrm{Ar}}$ for $G$ and its Pointrjagin dual $\widehat G$ coincide. That is, $$\#_{\mathrm{Ar}}(G)=\#_{\mathrm{Ar}}(\widehat G).$$
\end{axiom}

For our setting, we have the surjection $\mathbb A^n\twoheadrightarrow \mathbb A^n/(\mathbb A^n(g)+F^n)=H^1(F,g)$
with $H^1(F,g)$ compact, and hence an injection
$\widehat {H^1(F,g)}\hookrightarrow\widehat{\mathbb A^n}$ with $\widehat{H^1(F,g)}$ discrete.
Fix a test function $f$ on $\mathbb A^n$, then we obtain its Fourier transform $\widehat f$ using
the character $e^{-2\pi\sqrt {-1}\sum_{v\in S}\chi_v(*)}$
on $\widehat{\mathbb A^n}=\mathbb A^n.$ In particular, it makes sense to talk about 
$\widehat f|_{\widehat {H^1(F,g)}}$. Now, for the pairing $({\widehat {H^1(F,g)}}, { {H^1(F,g)}})$,
 using the character $e^{-2\pi\sqrt {-1}\sum_{v\in S}\chi_v(*)}$ again, we obtain the Fourier transform $\widetilde{\widehat f\ }$ of $\widehat f|_{\widehat {H^1(F,g)}}$ defined by
$$\widetilde{\widehat f\ }(\eta):=\sum_{\xi\in\widehat {H^1(F,g)}}e^{-2\pi\sqrt {-1}\sum_{v\in S}\chi_v(\xi_v\eta_v)}\widehat f(\xi),\qquad\forall\eta\in H^1(F,g).$$
Assume the test function $f$ on $\mathbb A^n$ satisfying 
 $$f|_{H^0(F,g)}\in L^2\big({H^0(F,g)}\big), \quad\widehat f|_{\widehat {H^1(F,g)}}
\in L^2\big(\widehat {H^1(F,g)}\big),\quad \widetilde{\widehat f\ }\in L^2\big(H^1(F,g)\big).$$

\noindent
 Then we can introduce the arithmetic counts for $H^0(F,g)$ and $H^1(F,g)$ with respect to the test function $f$. 
 That is to say, by defining 
$$\begin{aligned}\#_{\mathrm{Ar},f}\big(H^0(F,g)\big):=&\int_{H^0(F,g)}\big|f({\bf x})\big|^2d{\bf x}
=\sum_{{\bf x}\in H^0(F,g)}\big|f({\bf x})\big|^2;\\
\#_{\mathrm{Ar},f}\big(\widehat{H^1(F,g)}\big):=&\int_{\widehat{H^1(F,g)}}\big|\widehat f({\bf y})\big|^2d{\bf y}
=\sum_{{\bf y}\in \widehat {H^1(F,g)}}\big|\widehat f({\bf y})\big|^2;\\
\#_{\mathrm{Ar},f}\big(H^1(F,g)\big):=&\int_{H^1(F,g)}\big|\widetilde{\widehat f\ }({\bf z})\big|^2d{\bf z}.\end{aligned}$$
In particular, 
by the Plancherel formula,
we have Axiom 4: $$\#_{\mathrm{Ar},f}\big(\widehat{H^1(F,g)}\big)=
\#_{\mathrm{Ar},f}\big(H^1(F,g)\big).$$ 

\noindent
{\it Remark.} The above procedure can be applied to introduce counts for all discrete subgroups and  compact quotient groups of locally compact groups.

To define $h^1$, as above, we use the canonical test function $${{\bf e}_{\mathbb A}^{\frac{1}{2}}}:=\prod_{\frak p\in S_{\mathrm{fin}}}{\bf 1}_{\mathcal O_\frak p^n}\times\prod_{v\in S_\infty}e^{-\frac{1}{2}\pi N_v \|*\|_{\mathrm{can},v}}$$
on $\mathbb A^n$. (Recall that being the characteristic function, ${\bf 1}_{\mathcal O_\frak p^n}=
{{\bf 1}_{\mathcal O_\frak p^n}^{\frac{1}{2}}}$.)
We then have a natural count for $H^1(F,g)$ using ${\bf e}(*):=\big|\widetilde{\widehat{\bf e^{\frac{1}{2}}_{\mathbb A}(g*)}}\big|^2$:
$$\#_{\mathrm{Ar}}H^1(F,g):=\int_{H^1(F,g)}{\bf e}({\bf x})d{\bf x}.$$

\begin{mdef} Define the 
 1-st numerical  cohomology $h^1(F,g)$ for a matrix idele $g\in GL_n(\mathbb A)$ by 
 $$h^1(F,{\bf g}):=\log \Big(\#_{\mathrm{Ar}}(H^1(F,g))\Big)=\log\Big(\int_{H^1(F,{\bf g})}{\bf e}({\bf x})\,d\mu({\bf x})\Big).$$
 \end{mdef}
By a direct calculation ([Ta]), we have 
$$\widehat{{\bf e}_{\bf A}}=\prod_{\frak p\in S_{\mathrm{fin}}}{\bf 1}_{(\partial^{-1}_\frak p)^n}\times\prod_{v\in S_\infty}e^{-\pi N_v\|*\|_{\mathrm{can},v}}.$$ This, together with our arithmetical duality 
$H^1(F,g)\simeq\widehat{H^0(F,\kappa_F\cdot g^{-1})},$ implies that
$$\#_{\mathrm{Ar}}\widehat{H^1(F,g)}=\sum_{{\bf a}\in\widehat{H^1(F,G)}}\widehat{\bf e_{\mathbb A}}(g^{-1}{\bf a})=\#_{\mathrm{Ar}}H^0(F,\kappa_F\cdot g^{-1}),$$
because, in the definition of $\#_{\mathrm{Ar}}H^0(F,g)$, a shift by $g$ is used.
But, our counting system satisfies Axiom 4 (due to the Plancherel formula), 
$$\#_{\mathrm{Ar}}(H^1(F,g))=\#_{\mathrm{Ar}}\widehat{H^1(F,g)}.$$ 
This then proves (1) above, or the same, the numerical duality below.
Moreover, from definition of the group $H^0(F,g)$, particularly the local condition $g_v{\bf x}_v\in\mathcal O_p^n$, we see that
$$\#_{\mathrm{Ar}}H^0(F,g)=\sum_{{\bf a}\in F^n}\prod_{\frak p\in S_{\mathrm{fin}}}{\bf 1}_{\mathcal O_\frak p^n}(g_\frak p {\bf  a})\times\prod_{v\in S_\infty}e^{-\pi N_v\|g_v{\bf a}\|_{\mathrm{can},v}}=\sum_{{\bf a}\in F^n}{\bf e}_{\mathbb A}(g \cdot{\bf a}).$$ Similarly, 
$$\#_{\mathrm{Ar}}H^0(F,\kappa_F\cdot g^{-1})=\sum_{{\bf a}\in F^n}{\bf e}_{\mathbb A}(\kappa_F\cdot g^{-1}\cdot{\bf a})=\sum_{{\bf a}\in F^n}\widehat{\bf e}_{\mathbb A}( g^{-1}\cdot{\bf a}).$$ On the other hand, by Tate's Riemann-Roch ([Ta]), obtained from  the Poisson summation formula for adelic spaces $F^n\hookrightarrow \mathbb A^n$, 
we have $$\sum_{{\bf a}\in F^n}{\bf e}_{\mathbb A}(g\cdot{\bf a})=\frac{1}{\|\mathrm{det}g\|}\cdot
\sum_{{\bf a}\in F^n}\widehat{\bf e}_{\mathbb A}( g^{-1}\cdot{\bf a}).$$
That is to say,
$$\#_{\mathrm{Ar}}H^0(F,g)=\frac{1}{\|\mathrm{det}g\|}\cdot \#_{\mathrm{Ar}}H^0(F,\kappa_F\cdot g^{-1}).$$ Or equivalently, with the help of the numerical duality (1), we have
$$h^0(F,g)-h^1(F,g)=\chi(F,g).$$ 
This, together with the Arakelov Riemann-Roch ([La1,2]) $$\chi(F,g)=\mathrm{deg}(g)-\frac{n}{2}\cdot\log \Delta_F,$$ then implies (2) above, and hence completes the proof of the following fundamental:

\begin{mthm} For  a matrix idele $g\in GL_n(\mathbb A)$, we have

\noindent
(i) (Numerical Duality) $$\displaystyle{h^1(F,g)=h^0(F,\kappa_F\cdot  g^{-1});}$$

\noindent
(ii) (Arithmetic Riemann-Roch) 
$$\displaystyle{h^0(F,g)-h^1(F,g)=\mathrm{deg}(g)-\frac{n}{2}\cdot\log \Delta_F.}$$
\end{mthm}

\noindent
$\bullet$ {\bf Lattice Version}: In parallel, we have  a theory for  metrized  bundles. Let $\Lambda$ be an $\mathcal O_F$-lattice of rank $n$  in the metrized space
$\Big(\prod_{\sigma:\mathbb R}\mathbb R^n\times\prod_{\tau:\mathbb R}\mathbb C^n,\rho\Big).$  Define its {\it 0-th and 1-st cohomological groups} by setting
 $H^0(F,\Lambda)=\Lambda$  to be the lattice itself and $H^1(F,\Lambda):=\Big(\prod_{\sigma:\mathbb R}\mathbb R^n\times\prod_{\tau:\mathbb R}\mathbb C^n,\rho\Big)/\Lambda$ to be its compact quotient. Then, by Theorem 1, we have 

\begin{mmthm} (Arithmetic Duality) There is a canonical identitification
$$\widehat{H^1(F,\Lambda)}\simeq  H^0(F,\omega_F\otimes\Lambda^\vee),$$ where $\widehat{H^1(F,\Lambda)}$ is 
the Pontrjagin dual of $H^1(F,\Lambda)$, and $\omega_F$ is the canonical lattice of $F$ defined by the projective module associated with the inverse of the global differential module $\partial_F$ of $\mathcal O_F$ together with the standard matric.
\end{mmthm}
Recall that in algebraic geometry, for a divisor $D$ of an irreducible algebraic curve over the finite field $\mathbb F_q$, $$q^{\mathrm{dim}H^0(F,D)}=\#H^0(C,D)=\sum_{{\bf x}\in H^0(C,D)}{\bf 1}_{\bf x}.$$ Motivated by this, to introduce a numerical $h^0$, we start with  the natural weight function $e^{-\pi\sum_vN_v\|*\|_{\rho_v}}$. Thus being discrete, we obtain 
an arithmetic count for $H^0(F,\Lambda)$:
 $$\#_{\mathrm{Ar}}H^0(F,\Lambda):=\sum_{{\bf x}\in H^0(F,\Lambda)}e^{-\pi\sum_vN_v\|{\bf x}\|_{\rho_v}}.$$ This coincides with  $k^0(F,\Lambda)$ from arithmetic effectivity of [GS], and hence yields their well-known $h^0(F,\Lambda):=\log k^0(F,\Lambda)$. 
 
As for $h^1$, we introduce a counting function ${\bf e}(*)$ on $H^1(F,\Lambda)$ following the idelic discussion using Fourier transform
and $e^{-\pi\sum_vN_v\|*\|_v}$. This then gives the arithmetic count of $H^1(F,\Lambda)$ and the numerical $h^1(F,\Lambda)$ as follows:
$$\#_{\mathrm{Ar}}H^1(F,\Lambda):=\int_{{\bf x}\in H^1(F,\Lambda)}{\bf e}({\bf x})d\mu({\bf x})\ \ \mathrm{and}\ \ h^1(F,\Lambda):=\log \#_{\mathrm{Ar}}H^1(F,\Lambda).$$
Consequently, from the standard theory of Fourier analysis for lattices, we see that the topological duality and
the Plancherel theorem implies the numerical duality and the Poisson summation theorem (together with the numerical duality and the Arakelov Riemann-Roch) gives the Riemann-Roch:

\begin{mmthm}For an $\mathcal O_F$-lattice $\Lambda$ of rank $n$, we have

\noindent
(i) (Arithmetic Duality) $h^1(F,\Lambda)=  h^0(F,\omega_F\otimes\Lambda^\vee)$;

\noindent
(ii) (Riemann-Roch) $h^0(F,\Lambda)-h^1(F,\Lambda)=\mathrm{deg}(\Lambda)-\frac{n}{2}\cdot\log\Delta_F.$
\end{mmthm}
We end this subsection by drawing reader's attentions to [Bo], [Neu], [Mo], [Se] and [De].
 
\subsection{Ampleness and Vanishing Theorem}
Two reasons have made arithmetic vanishing theorem appeared difficult. First of all, in current theories, there is no individual arithmetical cohomology $h^i$ but rather a combined arithmetic Euler characteristic  
$\chi$; Secondly, even with the genuine arithmetic cohomology groups $H^i$ and hence $h^i$'s, it is impossible to have  zero 
 groups $H^0$ and $H^1$: After all, the genuine global $H^0$ (resp. $H^1$) in the case of number fields $F$ are discrete groups 
(resp.  compact groups) equal to (resp. dual to)  full rank lattices in $\big(\mathbb R^{r_1}\times\mathbb C^{r_2}\big)^n$ which can never become zero.

However, here, we, motivated by [GS] and [Gr], prove the following

\begin{mthm} (Vanishing Theorem) Let ${\bf a}\in GL_1(\mathbb A)$ be an idele of  $F$. Then $${\bf a}\ {is\ positive}\qquad\Leftrightarrow\qquad
\lim_{m\to \infty}h^1(F,  {\bf a}^m\cdot {\bf g})=0,\qquad\forall {\bf g}\in GL_n(\mathbb A).$$
\end{mthm}

By definition, an idele ${\bf a}=(a_\frak p;a_v)$ is called {\it positive} if $$\mathrm{deg}({\bf a}):=\sum_\frak p
\mathrm{ord}_\frak p(a_\frak p)\cdot N_\frak p\log p-\sum_vN_v\log|a_v|_v>0;$$
and ${\bf a}=(a_\frak p;a_v)$ is called {\it ample} if for each $g=(g_\frak p;g_v)
\in GL_n(\mathbb A)$,
the unit ball $B_0(1)$ centered at 0 in the metrized
space $\Big(\big(\mathbb R^{r_1}\times\mathbb C^{r_2}\big)^n,\rho(a_\infty^m\cdot g_\infty)\Big)$ contains a basis whose
 elements are positive at finite places of the lattice
$H^0(F,{\bf a}^m\cdot{\bf g})$, for sufficiently large $m\gg 0$ ([Zh]). 

\begin{mthm} (Criterions for Ampleness) Let ${\bf a}$ be an idele of a number field $F$. Then the following conditions are equivalent:

\noindent
(i)  ${\bf a}$ is ample;

\noindent
(ii)   ${\bf a}$ is positive;

\noindent
(iii)  $\displaystyle{\lim_{m\to \infty}h^1(F,  {\bf a}^m\cdot {\bf g})=0}$ for all ${\bf g}\in GL_n(\mathbb F).$
\end{mthm}

\noindent
{\it Proof.} (i) $\Rightarrow$ (ii). This is obvious.  

\noindent (ii) $\Rightarrow$ (i). This is essentially proved in ([Zh, Thm 2.2]).
A technical point here is that for a fixed $n$ and for all possible $m, {\bf a}$ and ${\bf g}$, 
$H^0(F,  {\bf a}^m\cdot {\bf g})\subset F^n$   are full rank lattices in metrized
 spaces based on (the vector space) $\big(\mathbb R^{r_1}\times\mathbb C^{r_2}\big)^n$.

Indeed, if ${\bf a}$ is positive,  by the Riemann-Roch, for sufficiently large $m$, $H^0(F,{\bf a}^m)$ contains a strictly effective section $l$, namely, a section $l$ satisfying $\|l\|<1$. Thus from what we just said, $H^0(F,  {\bf a}^m\cdot {\bf g})$ contains a full rank sublattice
generated by strictly effective sections $l^m\cdot e_1, l^m\cdot e_2,\dots, l^m\cdot e_N$, where  $N:=n\cdot [F:\mathbb Q]$, and 
$\{e_1,e_2,\dots, e_N\}$ is a fixed basis of $H^0(F, {\bf g})$. (Here, we view $l$ as an element of $F$ and $e_i$ as elements of $F^n$, and the product $\cdot$ is the scalar multiplication.)  This, together with Lemma 1.7 of [Zh], then shows that for sufficiently large $m$, $H^0(F,{\bf a}^m\cdot {\bf g})$ is generated by strictly effective sections. This establishes
the equivalence of  (i) and (ii).

\noindent
(iii) $\Rightarrow$ (ii). From (iii), by the duality, we have for $h^0(F,{\bf a}^{-m}\cdot {\bf g})\to 0$ as $m\to \infty$ 
for ${\bf g}\in GL_n({\mathbb A})$. 
But $$h^0(F,{\bf a}^{-m})=\log\Big(1+\sum_{{\bf x}\in H^0(F,{\bf a}^{-m}),{\bf x}\not=0}e^{-\pi\sum_vN_v\|{\bf x}\|_v}\Big).$$ Consequently, for any non-zero ${\bf x}(m)\in H^0(F,{\bf a}^{-m})$,
$$e^{-\pi\sum_vN_v\|{\bf x}(m)\|_v}\to 0,\quad\mathrm{or\ equivalently,}\ \|{\bf x}(m)\|_v\to \infty, \qquad \mathrm{as}\ m\to\infty.$$ Applying this to ${\bf x}(m):={\bf e}^m\in H^0(F,{\bf a}^{-m}),$  we conclude that $\|{\bf e}\|_v>1,\ \forall v$ for all non-zero sections ${\bf e}\in H^0(F,{\bf a}^{-1})$. This then further implies that  $\mathrm{deg}({\bf a})\geq 0$. 
 Indeed, if  $\mathrm{deg}({\bf a})< 0$ or the same $\mathrm{deg}({\bf a}^{-1})> 0$, 
 by applying the equivalence of (i) and (ii) to ${\bf a}^{-1}$, (replacing ${\bf a}$ with ${\bf a}^l$ for sufficiently large $l$ if necessary,) we conclude that $H^0(F,{\bf a}^{-1})$ consists of a $\mathbb Z$-basis $\{{\bf e}_1,{\bf e_2},\dots,{\bf e}_{n[F:\mathbb Q]}\}$ such that $\|{\bf e}_i\|_v< 1$,  a contradiction.

Moreover, we claim that $\mathrm{deg}({\bf a})\not=0$. Otherwise,  choose ${\bf g}_0$ such that 
$\chi(F,{\bf g}_0)=\mathrm{deg}({\bf g}_0)-\frac{n}{2}\log\Delta_F>0$. 
 Then by the Riemann-Roch,
$$h^0(F,{\bf a}^{-m}\cdot {\bf g}_0)-h^1(F,{\bf a}^{-m}\cdot {\bf g}_0)\Big(=\mathrm{deg}({\bf g}_0)-\frac{n}{2}\log\Delta_F\Big)=\chi(F,{\bf g}_0).$$
Consequently,  $\lim_{m\to \infty}h^0(F,{\bf a}^{-m}\cdot g)\geq\chi(F,{\bf g}_0)> 0$, a contradiction as well. This then implies that $\mathrm{deg}({\bf a})>0$, namely, (ii).
 
\noindent
(ii) $\Rightarrow$ (iii). By the numerical duality, it suffices to show that $$\lim_{m\to\infty}h^0(F,{\bf a}^{-m}\cdot{\bf g})=0\qquad\mathrm{ for\ any\ fixed}\ {\bf g}\in GL_n(\mathbb A).$$ Denote by $\lambda_1(-m)$ the first Minkowski successive minimum of the lattice $H^0(F,{\bf a}^{-m}\cdot{\bf g})\subset \Big(\mathbb R^{r_1}\times\mathbb C^{r_2}\Big)^n$.
Then, by Prop. 4.4 of [Gr], it suffices to show that $\lim_{m\to\infty}\lambda_1(-m)=+\infty$, since the $\mathbb Z$-rank of
the $\mathcal O_F$-lattice $H^0(F,{\bf a}^{-m}\cdot{\bf g})=H^0(F,{\bf a}^{-1})^{\otimes m}\otimes H^0(F,{\bf g})$ remains to be  $n\cdot[F:\mathbb Q]$. 

Suppose, otherwise, that there exists an increasing sequence $\{m_k\}$ of natural numbers such that $\lambda_1(-m_k)$ remains bounded. Replacing it with a subsequence if necessary, we may assume that $\lim_{k\to\infty}\lambda_1(-m_k)=\lambda_1$.

We claim that this is impossible. To see this, let us go as follows based on the stability to be introduced in the next section.
First of all, denote by $$\{0\}\subset\Lambda_1\subset\Lambda_2\subset\dots\subset\Lambda_s=H^0(F,{\bf g})$$  the Harder-Narasimhan filtration of the $\mathcal O_F$-lattice $H^0(F,{\bf g})$. Then the Harder-Narasimhan filtration of the lattice $H^0(F,{\bf a}^{-m}\cdot{\bf g})$ is given by 
$$\{0\}\subset H^0(F,{\bf a}^{-m})\cdot\Lambda_1\subset\dots\subset H^0(F,{\bf a}^{-m})\cdot\Lambda_s=H^0(F,{\bf a}^{-m}\cdot{\bf g}).$$
In particular,  $H^0(F,{\bf a}^{-m})\cdot\Lambda_1$ is,  by definition, a maximal semi-stable $\mathcal O_F$-sublattice of $H^0(F,{\bf a}^{-m}\cdot{\bf g})$ with biggest slope and highest rank, and  $$\lim_{m\to\infty}\mathrm{deg}\Big(H^0(F,{\bf a}^{-m})\cdot\Lambda_1\Big)=\lim_{m\to\infty}\Big(-m\cdot\mathrm{rank}(\Lambda_1)\cdot\mathrm{deg}({\bf a})+\mathrm{deg}(\Lambda_1)\Big)\to -\infty.$$ 

Secondly, fix a minimal vector ${\bf e}(-{m_k})\in H^0(F,{\bf a}^{-m_k}\cdot{\bf g})$ such that $\|{\bf e}(-{m_k})\|=\lambda_1(-m_k),$ and 
denote by  $\Lambda_1(-m_k)$ the rank one  $\mathcal O_F$-sublattices of $H^0(F,{\bf a}^{-m_k}\cdot{\bf g})$ generated by ${\bf e}(-{m_k})$.  Then, by the biggest slope and highest rank property of $H^0(F,{\bf a}^{-m})\cdot\Lambda_1$ in the
 the Harder-Narasimhan filtration of $H^0(F,{\bf a}^{-m_k}\cdot{\bf g})$,   we conclude that 
$$\mathrm{deg}\Big(\Lambda_1(-m_k)\Big)\leq \frac{\mathrm{deg}\big(H^0(F,{\bf a}^{-m_k})\cdot\Lambda_1\big)}{\mathrm{rank}\big(H^0(F,{\bf a}^{-m_k})\cdot\Lambda_1\big)}.$$ Consequently,
$$\lim_{k\to\infty}\mathrm{deg}\big(\Lambda_1(-m_k)\big)=-\infty.$$

Finally, applying Prop 7.1 of [Gr] to the rank one $\Lambda_1(-m_k)$, we see  that $$\lambda_1(-m_k)=\|{\bf e}(-{m_k})\|\geq [F:\mathbb Q]\cdot e^{-\frac{2}{[F:\mathbb Q]}\cdot{\mathrm{deg}\big(\Lambda_1(-m_k)\big)}}.$$ Hence, $\lambda_1(-m_k)$ are unbounded, a contradiction.
This completes the proof.
 
\section{Moduli Spaces of Semi-Stable Lattices}

\subsection{Stability}
Let $F$ be a number field with $\mathcal O_F$ the integer ring and $\Delta_F$ the absolute value of the discriminant.
By definition, an $\mathcal O_F$-lattice $\Lambda$ is called {\it semi-stable} if for any $\mathcal O_F$-sublattice $\Lambda'$ of 
$\Lambda$, we have $$\Big(\mathrm{Vol}(\Lambda)\Big)^{\mathrm{rank}(\Lambda')}\leq \Big(\mathrm{Vol}(\Lambda')\Big)^{\mathrm{rank}(\Lambda)},\quad\mathrm{or\ equivalently},\quad
\mu(\Lambda')\leq\mu(\Lambda),$$ where 
$\mu(\Lambda):=\frac{\mathrm{deg}(\Lambda)}{\mathrm{rank}_{\mathcal O_F}(\Lambda)}$; and
a size $n$ matrix idele  ${\bf g}\in\mathbb GL_n(\mathbb A)$ is called {\it semi-stable}  if the associated $\mathcal O_F$-lattice $\Lambda({\bf g})=H^0(F,{\bf g})$  is semi-stable. (The above equivalence is established using the Arakelov Riemann-Roch.)

Denote by $\mathcal M(n;d)=\mathcal M_F(n;d)$ (resp. $\mathcal M_{\mathrm{Ar},F}(n;d)$)
 the moduli space of semi-stable lattices of rank $n$ (resp. semi-stable matrix ideles of size $n$ and degree $d$). It is known that (see e.g., [Stu], [Gra] and [We3,4,6])

\noindent
(i) $\mathcal M_{\mathrm{Ar},F}(n;d)$ is a closed subset of 
 $GL_n(\mathbb A)$; 
  
 \noindent
 (ii) the fiber of the natural map $\mathcal M_{\mathrm{Ar}}(n;d)\to\mathcal M(n;d)$ is the compact  
 group $\prod_{\frak p\in S_{\mathrm{fin}}} GL_n(\mathcal O_{\frak p})
 \times\prod_{\sigma:\mathbb R}O_n(\mathbb R)\times\prod_{\tau:\mathbb C}U_n(\mathbb C)$;
 and 
 
 \noindent
 (iii) 
 $\mathcal M(n;d)$ and hence also  $\mathcal M_{\mathrm{Ar}}(n;d)$ are compact.
 
If $n=1$, the stability condition is automatic, and these compact spaces then coincide with
the well-known Arakelov  Picard groups
 $\mathrm{Pic}^d(F)$ and the idelic class group $\mathrm{Pic}^d_{\mathrm{Ar}}(F)$.
 
 \subsection{Effective Vanishing Theorem}
 
From now on, to avoid duplication, we present only the lattice version of the theory.
Since $h^0$ is a smooth function on the moduli spaces, it is quite natural for us to introduce the
{\it analytic stratifications} $\mathcal M(n,d)^{<(\mathrm{resp}.\leq,>,\geq,=)T}$ for a fixed real number $T$. For example,
 $$ \mathcal M(n,d)^{\leq T}:=\Big\{\Lambda\in
\mathcal M(n,d):h^0(F,\Lambda)\,\leq\,T\Big\}.$$
Since $\mathcal M(n,d)$ is compact, we have the following

\begin{prop} For fixed $n$ and $d$, there exist real numbers $M(n,d)$ and $m(n,d)$ such that\newline
\hskip 3.0cm $\mathcal M(n,d)^{< m}=\emptyset$ and $\mathcal M(n,d)^{> M}=\emptyset$\newline
if and only if $m,\,M\not\in [m(n,d),M(n,d)].$
\end{prop}

We call $m(n,d)$ and $M(n,d)$ the {\it minimal and maximal values} of $h^0$ on $\mathcal M(n,d)$ respectively, and $\delta(n,d):=M(n,d)-m(n,d)$ the {\it complicity} of $\mathcal M(n,d)$. 
 Since for any rank one lattice $\overline L$, $$\mathrm{det}( \Lambda^\vee\otimes \overline L)=
\mathrm{det}( \Lambda)^\vee\otimes (\overline L)^{\otimes{\mathrm{rank} \Lambda}},$$ by duality and the Riemann-Roch, we have the following

\begin{prop} (Duality) For extremal values of $h^0$ on the moduli spaces of semi-stable lattices, 
$$\begin{aligned}M\Big(n,n\log\Delta_F-d\Big)=&M(n,d)+\Big(\frac{n}{2}\log\Delta_F-d\Big),\\
m\Big(n,n\log\Delta_F-d\Big)=&m(n,d)+\Big(\frac{n}{2}\log\Delta_F-d\Big).\end{aligned} $$
\end{prop}

\begin{mthm}(Uniform Boundness) For any $\varepsilon>0$, these exists an effectively computable constant $d_F(n;\varepsilon)$ depending also 
on $F$ and $n$ such that, for all $d\geq d_F(n;\varepsilon)$
$$0<m_F(n;d)-\Big(d-\frac{n}{2}\log\Delta_F\Big)\leq M_F(n;d)-\Big(d-\frac{n}{2}\log\Delta_F\Big)<\varepsilon.$$
In particular, $$\lim_{d\to\infty}\delta_F(n;d)=0.$$
\end{mthm}

\noindent
{\it Proof.} Since the moduli spaces are compact, by the Riemann-Roch theorem, it suffices to establish the following:

\begin{mthm}(Effective Vanishing Theorem) Let $ \Lambda$ be a rank $n$ semi-stable $\mathcal O_F$-lattice.

\noindent
(i) If
$\displaystyle{\mathrm{deg}( \Lambda)\leq -[F:\mathbb Q]\cdot\frac{n\log n}{2}},$
we have
$$h^0( \Lambda)\leq \frac{3^{ n\cdot [F:\mathbb Q]}}{1-\frac{\log 3}{\pi}}\cdot \Big(e^{-\pi\cdot [F:\mathbb Q]}\Big)^{
 e^{-\frac{2\mathrm{deg}( \Lambda)}{ n\cdot [F:\mathbb Q]}}}.$$

\noindent
(ii) If
$\displaystyle{\mathrm{deg}( \Lambda)\geq [F:\mathbb Q]\cdot\frac{n\log n}{2}}+n\log\Delta_F,$
we have
$$h^1( \Lambda)\leq \frac{3^{ n\cdot [F:\mathbb Q]}}{1-\frac{\log 3}{\pi}}\cdot \Big(e^{-\pi\cdot [F:\mathbb Q]\cdot \Delta_F^{-\frac{2}{[F:\mathbb Q]}}
}\Big)^{ e^{\frac{2\mathrm{deg}( \Lambda)}{ n\cdot [F:\mathbb Q]}}}.$$
\end{mthm}
 
\noindent
{\it Remark.}
If $n=1$, (i) is proved in [Gr] (see also [GS]): the stability condition in rank one is automatic.

\noindent
{\it Proof.} By the numerical duality, it suffice to establish (i). For any non-zero vector ${\bf x}$ of $ \Lambda$, denote by $L({\bf x})=\mathcal O_F\cdot {\bf x}$, the $\mathcal O_F$-lattice generated by ${\bf x}$  in $ \Lambda$. Then by Lem 7.1 of [Gr], we know that $$\|{\bf x}\|^2\geq [F:\mathbb Q]\cdot e^{-\frac{2}{[F:\mathbb Q]}\cdot \mathrm{deg}(L({\bf x}))}.$$ But by the semi-stability of $ \Lambda$, we have
$$\mathrm{deg}(L({\bf x}))\leq \frac{\mathrm{deg}( \Lambda)}{n}\leq -[F:\mathbb Q]\cdot\frac{\log n}{2}.$$
Consequently, $\|{\bf x}\|^2\geq [F:\mathbb Q]\cdot n$. Therefore, Prop 4.4 of [Gr] can be applied to the $\mathbb Z$-lattice $ \Lambda$.
Note that from above, (by semi-stability,) $$\|{\bf x}\|^2\geq [F:\mathbb Q]\cdot e^{-\frac{2\mathrm{deg}( \Lambda)}{ n\cdot [F:\mathbb Q]}}
.$$ This implies $$\lambda_1^2\geq  [F:\mathbb Q]\cdot e^{-\frac{2\mathrm{deg}( \Lambda)}{ n\cdot [F:\mathbb Q]}}
.$$ Consequently, by Prop. 4.4 of [Gr], we have $$\#_{\mathrm{Ar}}H^0(F, \Lambda)
\leq 1+ \frac{3^{ n\cdot [F:\mathbb Q]}}{1-\frac{\log 3}{\pi}}
\cdot e^{-\pi\cdot [F:\mathbb Q]\cdot e^{-\frac{2\mathrm{deg}( \Lambda)}{ n\cdot [F:\mathbb Q]}}}.$$ Clearly,
$h^0(F, \Lambda)=\log \#_{\mathrm{Ar}}H^0(F, \Lambda)\leq  \#_{\mathrm{Ar}}H^0(F, \Lambda)-1.$ This completes the proof.

\subsection{Arithmetic Stratifications}

Besides the above analytic stratifications, we can also introduce  arithmetic stratifications  for these 
moduli spaces. To explain this, let us assume that $F_0\subset F$ is a subfield of $F$. Then an $\mathcal O_F$-lattice $\Lambda$ may be naturally viewed as an $\mathcal O_{F_0}$-lattice of rank $\mathrm{rank}_F(\Lambda)\cdot[F:F_0]$, which we denote by $\mathrm{Res}_{F_0}^F(\Lambda)$. Fix a convex polygon $g$. Then as $\mathcal O_{F_0}$-lattices,
$\mathrm{Res}_{F_0}^F(\Lambda)$ admits a natural Harder-Narasimhan type filtration ([Stu], [Gra], [We3,4]). Denote its associated canonical polygon by $\overline g_{F_0}(\Lambda)$, and introduce {\it arithmetical stratifications} by  $$\mathcal M_F(n,d)^{\leq_{F_0} g}:=\Big\{\Lambda\in\mathcal M_F(n,d)\,|\,\overline g_{F_0}(\Lambda)\leq g\Big\}.$$ This is a much more refined stability: Unlike the $F$-stability, for elements in  $\mathrm{Pic}^d(F)$, 
$F_0$-stability is far from being trivial. 
In fact, we expect these $F_0$-level canonical semi-stable filtration and arithmetical stratifications play key roles in the studies of moduli spaces of semi-stable bundles, say, in finding the arithmetical analogues for results in classical algebraic geometry related to special divisors on curves ([GS], [Gr], [Fr]).
  
\section{Non-Abelian Zeta Functions}

Let $F$ be a number field with $\Delta_F$ the absolute value of the discriminant of $F$,
and denote $\mathcal M_F(n):=\cup_{d\in\mathbb R}\mathcal M_F(n,d)$  the moduli space of rank $n$ semi-stable $\mathcal O_F$-lattices. The natural Tamagawa measure on $GL_n(\mathbb A)$ induces a natural measure on $\mathcal M_F(n)$ which we write as $d\mu$.

\begin{mdef} The rank $n$ non-abelain zeta function $\widehat{\zeta}_F(s)$ of $F$ is the integration
$$\widehat{\zeta}_{F,n}(s):=\Big(\Delta_F^{\frac{n}{2}}\Big)^s\cdot\int_{\Lambda\in\mathcal M_F(n)}\Big(e^{h^0(F,\Lambda)}-1\Big)\cdot\Big(e^{-s}\Big)^{\mathrm{deg}(\Lambda)}d\mu(\Lambda),\ \ \mathrm{Re}(s)>1.$$
\end{mdef}
\begin{mthm}  (0) {\it $\widehat{\zeta}_{F,1}(s)\buildrel\cdot\over=\widehat\zeta_F(s)$ the completed Dedekind zeta function;}

\noindent
 (1) ({\bf Meromorphic Continuation})
 {\it $\widehat{\zeta}_{F,n}(s)$  is well-defined when $\mathrm{Re}(s)>1$ and admits a
 meromorphic continuation, denoted also by 
$\widehat{\zeta}_{F,n}(s)$, to the whole complex $s$-plane;}

\noindent
(2) ({\bf Functional Equation}) $\widehat{\zeta}_{F,n}(1-s)=\widehat{\zeta}_{F,n}(s)$;

\noindent
(3) ({\bf Singularities\,\&\,Residues})  {\it $\widehat{\zeta}_{F,n}(s)$ has two
 singularities, all simple poles, at $s=0,1$,
 with the residues $\pm\mathrm{Vol}\big({\mathcal M}_{F,n}[1]\big)$,
where ${\mathcal M}_{F,n}[1]$ denotes 
the moduli space of rank $n$ semi-stable 
$\mathcal O_F$-lattices of volume one.}
\end{mthm}

\noindent
{\it Proof.} (0) is essentially due to Iwasawa ([Iw]) and Tate ([Ta]). Write the volume $T$ (resp. $\geq T$, resp. $\leq T$) part of the moduli space $\mathcal M_F(n)$ as
$\mathcal M_{F,n}[T]$ (resp. $\mathcal M_{F,n}[\geq T]$, resp. $\mathcal M_{F,n}[\leq T]$), then we have natural decompositions $$\mathcal M_F(n)=\mathcal M_{F,n}[\leq T]\cup \mathcal M_{F,n}[\geq T]=\cup_{T> 0}\mathcal M_{F,n}[T]$$ and $$d\mu=\frac{dT}{T}\cdot d\mu_T,$$ where $d\mu_T$ denote the natural induced volume form on $\mathcal M_{F,n}[T]$. 
Then $$\begin{aligned}\widehat{\zeta}_{F,n}(s)=&\int_{\Lambda\in\mathcal M_F(n)}\Big(e^{h^0(F,\Lambda)}-1\Big)\cdot\mathrm{Vol}(\Lambda)^s\cdot d\mu(\Lambda)\\
=&I(s)+A(s)-\alpha(s)\end{aligned}$$ where
$$\begin{aligned}
I(s):=&\int_{\Lambda\in\mathcal M_{F,n}[\geq 1]}\Big(e^{h^0(F,\Lambda)}-1\Big)\cdot\mathrm{Vol}(\Lambda)^s\cdot d\mu(\Lambda)\\
A(s):=&\int_{\Lambda\in\mathcal M_{F,n}[\leq 1]}e^{h^0(F,\Lambda)}\cdot\mathrm{Vol}(\Lambda)^s\cdot d\mu(\Lambda)\\
\alpha(s):=&\int_{\Lambda\in\mathcal M_{F,n}[\leq 1]}\mathrm{Vol}(\Lambda)^s\cdot d\mu(\Lambda).\end{aligned}$$

\noindent
{\bf I(s)}: By the effective vanishing theorem, $I(s)$ is  holomorphic  over the whole complex $s$-plane, since it is the  integration of 
$T$-exponentially decay function $\Big(e^{h^0(F,\Lambda)}-1\Big)\cdot\mathrm{Vol}(\Lambda)^s$  
over the space $\mathcal M_{F,n}[\geq 1]=\cup_{T\geq 1}\mathcal M_{F,n}[T]$;
\vskip 0.20cm
\noindent
{\bf A(s)}: Note that if $\Lambda$ is semi-stable, so is $\kappa_F\otimes \Lambda^\vee$. Consequently,
$\Lambda\mapsto\kappa_F\otimes \Lambda^\vee$ defines a natural involution on $\mathcal M_F(n)$, and interchanges 
$\mathcal M_{F,n}[\geq 1]$  and $\mathcal M_{F,n}[\geq 1]$ by the duality and the Riemann-Roch.
Thus
$$\begin{aligned}A(s):=&\int_{\Lambda\in\mathcal M_{F,n}[\leq 1]}e^{h^0(F,\Lambda)}\cdot e^{-s\cdot \chi(F,\Lambda)}\cdot d\mu(\Lambda)\\
=
&\int_{\Lambda\in\mathcal M_{F,n}[\geq 1]}e^{h^0(F,\kappa_F\otimes\Lambda^\vee)}\cdot
\cdot e^{-s\cdot \chi(F,\kappa_F\otimes\Lambda^\vee)} d\mu(\Lambda)\\
=
&\int_{\Lambda\in\mathcal M_{F,n}[\geq 1]}e^{h^1(F,\Lambda)}\cdot e^{s\cdot \chi(F,\Lambda)}
 d\mu(\Lambda)\\
 =&\int_{\Lambda\in\mathcal M_{F,n}[\geq 1]}e^{h^0(F,\Lambda)}\cdot e^{(s-1)\cdot \chi(F,\Lambda)}
 d\mu(\Lambda)\\
 =&\int_{\Lambda\in\mathcal M_{F,n}[\geq 1]}e^{h^0(F,\Lambda)}\cdot \mathrm{Vol}(\Lambda)^{1-s} d\mu(\Lambda)\\
 =&I(1-s)+\beta(s)\end{aligned}$$
 where $$\beta(s):=\int_{\Lambda\in\mathcal M_{F,n}[\geq 1]}\mathrm{Vol}(\Lambda)^{1-s}\cdot d\mu(\Lambda).$$ Since $I(1-s)$ is a holomorphic function on $s$,
it suffices to understand $\beta(s)$.
\vskip 0.20cm
\noindent
${\bf \alpha(s)}$ {\bf and} ${\bf \beta (s)}$: By definition,
$$\begin{aligned}\alpha(s)=&\int_{\Lambda\in\mathcal M_{F,n}[\leq 1]}\mathrm{Vol}(\Lambda)^s\cdot d\mu(\Lambda)\\
=&\int_0^1 \frac{dT}{T}\int_{\Lambda\in \mathcal M_{F,n}[T]}\mathrm{Vol}(\Lambda)^s\cdot d_T\mu(\Lambda)\\
=&\int_0^1 \frac{dT}{T}\int_{\Lambda\in \mathcal M_{F,n}[T]}T^s\cdot d_T\mu(\Lambda)\\
=&\int_0^1 T^s\frac{dT}{T}\cdot\int_{\Lambda\in \mathcal M_{F,n}[T]} d_T\mu(\Lambda)\end{aligned}$$
Note that there is a natural isomorphism between $\mathcal M_{F,n}[1]\to \mathcal M_{F,n}[T]$, say by sending $\Lambda\mapsto \Lambda\otimes\Lambda_0$ for a fixed rank one lattice $\Lambda_0$ of degree $\frac{\log T}{n}$.
$$\begin{aligned}\alpha(s)=&\int_0^1 T^s\frac{dT}{T}\cdot\mathrm{Vol}\Big(\mathcal M_{F,n}[1]\Big)\\
=&\mathrm{Vol}\Big(\mathcal M_{F,n}[1]\Big)\cdot\frac{1}{s}.\end{aligned}
 $$ Similarly,
 
 $$\begin{aligned}\beta(s)=&\int_{\Lambda\in\mathcal M_{F,n}[\geq 1]}\mathrm{Vol}(\Lambda)^{1-s}\cdot d\mu(\Lambda)\\
 =&-\mathrm{Vol}\Big(\mathcal M_{F,n}[1]\Big)\cdot\frac{1}{1-s}.\end{aligned}
 $$  All in all, the up shot is
 $$\begin{aligned}\widehat\zeta_{F,r}(s)=&I(s)+I(1-s)-\alpha(s)+\beta(s)\\
 =&I(s)+I(1-s)
 +\mathrm{Vol}\Big(\mathcal M_{F,n}[1]\Big)\cdot\Big(
 \frac{1}{s-1}-\frac{1}{s}\Big),\end{aligned}$$
 with $$I(s)=\int_{\Lambda\in\mathcal M_{F,n}[\geq 1]}\Big(e^{h^0(F,\Lambda)}-1\Big)\cdot\mathrm{Vol}(\Lambda)^s\cdot d\mu(\Lambda)$$ a holomorphic function of $s$.
 This then proves the theorem.
\vskip 0.20cm 
High rank zetas, the natural non-abelain counterparts of Dedekind zeta functions, are expect to play a central role in the study of non-coummutative arithmetic aspects of number fields.
\vskip 0.30cm
We conclude this paper with the following comments:
While our method  is a continuation of the classical one due to Chevalley, Weil ([W], see also [Se]), Iwasawa ([Iw]) and Tate ([Ta]), these genuine global arithmetic cohomologies and non-abelian zetas  expose  new structures for number fields ([We 1-6]). In particular, the study of the so-called abelian parts of our non-abelian zetas is now actively carrying on. For details, please refer to [H], [Ki], [KKS], [Ko], [LS], [Su1,2], [SW] and [We 2-6].\footnote{
I would like to thank C. Deninger, I. Fesenko, H. Hida, H. Kim, I. Nakamura, K. Ueno, H. Yoshida and D. Zagier for their  interests and encouragements, and  G. van der Geer for introducing me
his joint work with R. Schoof.}

\vfill
\eject
\vskip 1.0cm
\centerline {\bf\large{ REFERENCES}}
\vskip 0.40cm 
\noindent
[Bo] A. Borisov, Convolution structures and arithmetic cohomology, Compositio Math., 136 (2003), no. 3, 237-254.
\vskip 0.20cm  
\noindent
[De]  C. Deninger, Some analogies between number theory
and dynamical systems on foliated spaces. {\it Proceedings of the International Congress of
Mathematicians}, Berlin (1998)  Doc. Math.  1998,  Extra Vol. I, 163--186
 \vskip 0.20cm  
\noindent 
[Fr] P. Francini, The size function $h^0$ for quadratic number fields. 
Journal de Theorie des Nombres de Bordeaux, tome 13 , no 1 (2001), p. 125-135.
 \vskip 0.20cm  
\noindent
[GS] G. van der Geer \& R. Schoof, Effectivity of Arakelov
Divisors and the Theta Divisor of a Number Field, Sel. Math., New ser.
{\bf 6} (2000), 377-398  
 \vskip 0.20cm  
\noindent
[Gra] D. R. Grayson, Reduction theory using semistability. 
Comment. Math. Helv.  {\bf 59}  (1984),  no. 4, 600--634.
 \vskip 0.20cm  
\noindent
[Gr] R.P. Groenewegen, An arithmetic analogue of Clifford's Theorem. 
Journal de Theorie des Nombres de Bordeaux 13 (2001) 143-156.
 \vskip 0.20cm  
\noindent
[H] T. Hayashi, Computation of Weng's rank 2 zeta function over an algebraic number field,
J. Number Theory 125 (2007), no.2, 473-527
\vskip 0.20cm  
\noindent
[Iw] K. Iwasawa, Letter to Dieudonn\'e, April 8, 1952, 
 {\it Zeta Functions in Geometry}, 
Advanced Studies in Pure Math. {\bf 21} (1992), 445-450
 \vskip 0.20cm  
\noindent
[Ki] H. Ki, On the zeros of Weng's zeta functions,
Int Math Res Notices (2010) 2010, 2367-2393
 \vskip 0.20cm  
\noindent
[KKS] H. Ki, Y. Komori \& M. Suzuki, On the zeros of Weng zeta functions for Chevalley groups, preprint, 2010, arXiv:1011.4583
 \vskip 0.20cm  
\noindent
[Ko] Y. Komori, Functional equations for Weng's zeta functions for $(G,P)/\mathbb{Q}$, preprint, 2010, arXiv:1011.4582
 \vskip 0.20cm  
\noindent
[LS] J.C. Lagarias \& M. Suzuki, The Riemann hypothesis for certain integrals of Eisenstein series, J. Number Theory 118 (2006), no. 1, 98--122.
 \vskip 0.20cm  
\noindent
[La] S. Lang, {\it Algebraic Number Theory}, 
Springer-Verlag, 1986
 \vskip 0.20cm  
\noindent
[Mo] C. Moreno, {\it  Algebraic curves over finite fields.}
Cambridge Tracts in Mathematics, {\bf 97}, Cambridge University Press, 1991
 \vskip 0.20cm  
\noindent
[Neu] J. Neukirch, {\it Algebraic Number Theory}, Grundlehren der
Math. Wissenschaften, Vol. {\bf 322}, Springer-Verlag, 1999
 \vskip 0.20cm  
\noindent
[Se] J.-P. Serre, {\it Algebraic Groups and Class Fields}, GTM
{\bf 117}, Springer-Verlag (1988)
 \vskip 0.20cm  
\noindent
[Stu1] U. Stuhler, Eine Bemerkung zur Reduktionstheorie
quadratischer Formen.  Arch. Math.  {\bf 27} (1976), no. 6,
604--610 
 \vskip 0.20cm  
\noindent
[Stu2] U. Stuhler, Eine Bemerkung zur Reduktionstheorie
quadratischer Formen.   II.  Arch. Math.   {\bf 28}  (1977), no. 6, 611--619. 
 \vskip 0.20cm  
\noindent
[Su1] M. Suzuki, The Riemann hypothesis for Weng's zeta function of Sp(4) over Q, with an appendix by L. Weng, J. Number Theory 129 (2009), no. 3, 551-579. 
 \vskip 0.20cm  
\noindent
[Su2] M. Suzuki, A proof of the Riemann hypothesis for the Weng zeta function of rank 3 for the rationals. in {\it The Conference on L-Functions}, 175--199, World Sci. Publ., Hackensack, NJ, 2007. 
 \vskip 0.20cm  
\noindent
[SW] M. Suzuki, \& L. Weng, Zeta functions for $G_2$ and their zeros, {\it International Mathematics Research Notice}, 2009, 241-290
 \vskip 0.20cm  
\noindent
[Ta] J. Tate, Fourier analysis in number fields and Hecke's
zeta functions, Thesis, Princeton University, 1950 
 \vskip 0.20cm 
\noindent
[W] A. Weil, {\it Basic Number Theory},  Springer Varlag, 1973
\vskip 0.20cm 
\noindent
[We1] L. Weng, Non-abelian zeta functions for function fields. {\it Amer. J. Math}. 127 (2005), no. 5, 973--1017.
 \vskip 0.20cm 
\noindent 
[We2] L. Weng, A rank two zeta and its zeros. {\it J. Ramanujan Math. Soc.} 21 (2006), no. 3, 205--266
 \vskip 0.20cm 
\noindent 
[We3] L. Weng, Geometric arithmetic: a program, in {\it Arithmetic geometry and number theory,} 211--400, Ser. Number Theory Appl., 1, 
World Sci. Publ., Hackensack, NJ, 2006.
 \vskip 0.20cm 
\noindent 
[We4] L. Weng, A geometric approach to $L$-functions. {\it The Conference on $L$-Functions,} 219--370, World Sci. Publ., Hackensack, NJ, 2007.
 \vskip 0.20cm 
\noindent 
[We5] L. Weng, Symmetries and the Riemann Hypothesis, {\it Advanced Studies in Pure Mathematics 58}, Japan Math. Soc., 173-224, 2010
 \vskip 0.20cm 
\noindent 
[We6] L. Weng, Stability and Arithmetic,  {\it Advanced Studies in Pure Mathematics 58}, Japan Math. Soc., 225-360, 2010
 \vskip 0.20cm 
\noindent
[Zh] S. Zhang, Positive line bundles on arithmetic surfaces, Ann. Math., 136 (1992),
569-587.
\vskip 1.50cm
Lin WENG




Graduate School of Mathematics

Kyushu University

Fukuoka, 819-0395 

JAPAN

E-Mail: weng@math.kyushu-u.ac.jp

\end{document}